\documentclass[11pt,dvipsnames]{article}
\usepackage[utf8]{inputenc}
\usepackage[T1]{fontenc}
\usepackage{comment}
\usepackage{authblk}
\usepackage{relsize}
\usepackage{changepage}
\usepackage{mathtools}
\usepackage{booktabs}
\usepackage{multirow}
\usepackage{array}
\usepackage{enumerate}
\usepackage{enumitem}
\usepackage{mathtools}
\usepackage[english]{babel}
\usepackage{amsmath}
\usepackage{mathrsfs}
\usepackage{amsthm}
\usepackage{amsfonts}
\usepackage{url}
\usepackage{amssymb}
\usepackage{graphicx}
\usepackage{tabularx}
\usepackage[numbers]{natbib}
\usepackage{tikz}
\usetikzlibrary{positioning}
\usepackage{color}
\theoremstyle{plain}
\newtheorem{theorem}{Theorem}[section]
\newtheorem{corollary}[theorem]{Corollary} 
\newtheorem{defn}[theorem]{Definition}
\newtheorem{proposition}[theorem]{Proposition}
\newtheorem{lemma}[theorem]{Lemma}

\makeatletter 
\newcommand{\vast}{\bBigg@{4}} 
\newcommand{\Vast}{\bBigg@{5}} 
\makeatother 
\definecolor{bulgarianrose}{rgb}{0.28, 0.02, 0.03} 
\definecolor{gray}{rgb}{0.5, 0.5, 0.5} 
 
\theoremstyle{definition}

\newtheorem*{conjecture}{Conjecture}

\theoremstyle{remark} 
 
\usepackage[left=2cm,right=2cm,top=2cm,bottom=2cm]{geometry} 
\usepackage{amssymb} 
\makeatletter 
\def\namedlabel#1#2{\begingroup
    #2%
    \def\@currentlabel{#2}%
    \phantomsection\label{#1}\endgroup
} 
\usepackage[normalem]{ulem} 
\usepackage{dsfont} 
\usepackage{accents} 
\usepackage{verbatim} 
\usepackage{ulem} 
\usepackage{pgf,tikz,pgfplots} 
\pgfplotsset{compat = 1.16} 
\usepackage[all]{xy} 
 
\usepackage{stackengine} 
\stackMath 
\newcommand\tsup[2][2]{%
 \def\useanchorwidth{T}%
  \ifnum#1>1%
    \stackon[-.5pt]{\tsup[\numexpr#1-1\relax]{#2}}{\scriptscriptstyle\sim}%
  \else%
    \stackon[.5pt]{#2}{\scriptscriptstyle\sim}%
  \fi%
}

\usepackage{enumerate} 
 
\usepackage{amsmath} 
\usepackage{amsfonts} 
\usepackage{amsthm} 
\usepackage{amssymb} 
\usepackage{array} 
\usepackage{booktabs} 
\usepackage{changepage} 
\usepackage[english]{babel} 
\usepackage{enumerate} 
\usepackage{epsfig} 
\usepackage{float} 
\usepackage{framed} 
\usepackage{gensymb} 
\usepackage{graphicx} 
\usepackage[colorlinks=true, allcolors=.]{hyperref}
\usepackage{indentfirst} 
\usepackage{longtable} 
\usepackage{mathtools} 
\usepackage{multirow} 
\usepackage[normalem]{ulem} 
\usepackage{setspace} 
\usepackage{subcaption} 
\usepackage{textcomp} 
\usepackage[utf8]{inputenc} 
\usepackage{verbatim} 
\usepackage{enumerate} 
\usepackage{xlop} 
\usepackage{listings} 
\usepackage[none]{hyphenat}

\newcommand{\Mod}[1]{\ (\mathrm{mod}\ #1)}

\def\and{%
  \end{tabular}%
  \hskip 1em \@plus.17fil%
  \begin{tabular}[t]{c}}%
  
\title{\scshape On the maximal size of $(a,b)$-town$\Mod k$ families}  
 
\author{Nikola Veselinov, Miroslav Marinov\footnote{nikola.veselinov.veselinov@gmail.com, m.marinov1617@gmail.com} }
\date{} 
\begin{document}

 \maketitle

\begin{abstract}  
    A family $\mathcal{F}\subseteq\mathscr{P}(n)$ is an $(a,b)$-town$\Mod k$ if all sets in it have cardinality $a\Mod k$ and all pairwise intersections in it have cardinality $b\Mod k$. For $k=2$ the maximal size of such a family is known for each $a,b$, while for $k=3$ only $b-a\equiv 2 \Mod 3$ is fully understood. We provide a bound for $k=3$ when $b-a\equiv 1 \Mod 3$ and $n\equiv 2 \Mod 3$, which turns out to be tight for infinitely many such $n$. We also give sufficient conditions on the parameters $a,b,k,n$, which result in a better bound than the one from general settings by Ray-Chaudhuri--Wilson, in particular showing that this bound occurs infinitely often in a sense where all of $a,b,n$ can vary for a fixed $k$.
\end{abstract}

\section{Introduction}
\emph{Oddtown} in extremal set theory is a problem from the origins of algebraic combinatorics where elegant accessible proofs rely on linear algebra rather than pure combinatorial arguments. Throughout we denote by $[n]$ the set $\{1, 2, \ldots, n\}$ for a positive integer $n$ (conventionally, $[0]=\emptyset$), and by $\mathscr{P}(n)$ the power set of $[n]$. The original formulation of Oddtown concerns the maximal size of a family $\mathcal{F}\subseteq\mathscr{P}(n)$ whose sets have odd cardinality and whose intersections of any two distinct sets have even cardinality. This maximum was conjectured to be $n$ by Erd\H{o}s and proven independently by Berlekamp \cite{E-Berlekamp} and Graver \cite{J-Graver}. Their methods are applicable in more general settings -- the only importance of the integer $2$ is that it is a prime. More recently, a counting proof of the Oddtown problem has been established by Petrov \cite{Fedor-Petrov}, although the methods used are based on symmetric difference arguments and so are limited to mod $2$. 

A different generalization, also attracting considerable attention, is where individual set cardinalities are only required to be non-multiples of $k$ and pairwise intersections to be multiples of $k$. The extremal size of these families is proven to be $n$ for all prime power moduli, but is yet open for all others. The best known upper bound \cite{L-Babai-Frankl-Book} for $\operatorname{mod}\, 6$ is $2n-\log_2n$.

We consider the following generalization.

\begin{defn}
\label{defn:a-b-town-mod-k}
Let $a,b,k$ be integers with $k\geq 2$ and $0\leq a,b \leq k-1$. A family of sets $\mathcal{F}\subseteq\mathscr{P}(n)$ is an $(a,b)$-town$\Mod k$ if for any distinct $A, B\in\mathcal{F}$ we have $|A|\equiv a\Mod k$ and $|A\cap B|\equiv b\Mod k$.
\end{defn}

The classical Oddtown is $(1,0)$-town$\Mod 2$. The other three cases of $(a,b)$-town$\Mod 2$ are also well known -- for $(0,1)$ the exact result is $n-1+(n \Mod 2)$, while for $(0,0)$ (\emph{Eventown}) and $(1,1)$ the extremal sizes are $2^{\lfloor\frac{n}{2}\rfloor}$ and $2^{\lfloor\frac{n-1}{2}\rfloor}$, respectively \cite[Chapter 1 and 2]{L-Babai-Frankl-Book}.

We remark that if $\mathcal{F}$ is an $(a,b)$-town$\Mod k$, then it is also an $(a,b)$-town$\Mod d$ for any divisor $d$ of $k$. All results in this paper can be adapted to mod $k$ if we work with a single prime divisor of $k$, so throughout we shall work only with prime modulus $p$.

As it turns out, the maximal size of an $(a,b)$-town$\Mod p$ is exponential if and only if $a=b$. When $a\neq b$, a classical result by Ray-Chaudhuri--Wilson implies a linear upper bound in our setting.

\begin{theorem}\emph{(Modular RW \cite[Theorem 5.37]{L-Babai-Frankl-Book})}
\label{thm:modular-rw-theorem}
    Let $n$ be a positive integer, $p$ be a prime and $L$ be a set of $s\leq p-1$ integers. Let $t \geq 0$ be an integer with $t\notin L\Mod p$ and $s+t\leq n$. Let $\mathcal{F}\subseteq\mathscr{P}(n)$ be a family of sets such that $|E|\equiv t\Mod p$ and $|E\cap F|\in L\Mod p$ for any distinct $E,F\in\mathcal{F}$. Then $\displaystyle |\mathcal{F}|\leq {n\choose s}$.
\end{theorem}

For $s=1$, $t = a$, $L = \{b\}$ this implies $|\mathcal{F}| \leq n$ for any $(a,b)$-town $\Mod p$ if $p$ does not divide $a-b$.

\smallskip

Results concerning $(a,b)$-town are currently limited even for modulo $3$. Three of the nine cases have been solved \cite{T-ITYM-Germany-South}, namely, $(1,0)$, $(2,1)$ and $(0,2)$. For $(0,1), (1,2), (2,0)$ we have the linear bound from Theorem~\ref{thm:modular-rw-theorem}, but we are also able to improve it when $n\equiv 2 \Mod 3$ and partially for $n\equiv 0 \Mod 3$.

\begin{proposition}
\label{prop:mod-3-n-1}
    Suppose $\mathcal{F}$ is an $(a,b)$-town mod $3$ family of sets in $[n]$ where $b-a \equiv 1 \Mod 3$ and one of the following holds$:$
    \begin{itemize}
        \item $n\equiv 2 \Mod 3;$
        \item $n\equiv 0 \Mod 3$ and $a=0$, $b=1$.
    \end{itemize}
    Then $|\mathcal{F}| \leq n-1$.
\end{proposition}

This bound turns out to be essentially tight -- Lahtonen \cite{Jyrki-Math-Stack} has given a $(2,0)$-town family of size $n-1$ for infinitely many $n$, namely when $n\equiv 8\Mod{12}$ and $n-1$ is power of a prime, as well as a $(1,2)$-town of size $n-1$ when $n\equiv 7 \Mod {12}$ is a power of a prime. Furthermore, by taking the complement of each set in his examples, we obtain $(0,1)$-towns of size $n-1$ for infinitely many $n\equiv 1,2 \Mod 3$. 

Regarding $(0,0) \Mod p$, a known bound by Frankl and Odlyzko \cite{FRANKL-Odlyzko} is $2^{\left\lfloor\frac{n}{2}\right\rfloor}$. We are able to give a modified proof to show a similar one holds for $(m,m) \Mod p$.

\begin{proposition}
\label{prop:type-1-2-n-2-bound-intro}
    Let $p$ be a prime and $m$ be an integer with $1 \leq m \leq p-1$. For any $(m,m)$-town$\Mod p$ family $\mathcal{F} \subseteq \mathscr{P}(n)$ we have $|\mathcal{F}|\leq2^{\left\lfloor\frac{n+1}{2}\right\rfloor}$.
\end{proposition}

Hence so far the bounds on the maximal size for $(a,b)$-town$\Mod 3$ are as in the table below.

{\renewcommand{\arraystretch}{1.2} 

\begin{table}[ht]
\centering
\begin{tabular}{|c|c|c|c|}
\hline
$(a,b)$ &
$n \equiv 0\Mod 3$ &
$n \equiv 1\Mod 3$ &
$n \equiv 2\Mod 3$ \\
\hline
\hline
\multirow{2}{*}{$(0,0)$} &
Lower: $24^{\left\lfloor\frac{n}{12}\right\rfloor}$ & Lower: $24^{\left\lfloor\frac{n}{12}\right\rfloor}$ & Lower: $24^{\left\lfloor\frac{n}{12}\right\rfloor}$ \\
& Upper: $2^{\left\lfloor\frac{n}{2}\right\rfloor}$ & Upper: $2^{\left\lfloor\frac{n}{2}\right\rfloor}$ & Upper: $2^{\left\lfloor\frac{n}{2}\right\rfloor}$ \\
\hline
$(m,m)$ &
Lower: $24^{\left\lfloor\frac{n-m}{12}\right\rfloor}$ & Lower: $24^{\left\lfloor\frac{n-m}{12}\right\rfloor}$ & Lower: $24^{\left\lfloor\frac{n-m}{12}\right\rfloor}$ \\
$m \in \{1,2\}$ & Upper: $2^{\left\lfloor\frac{n+1}{2}\right\rfloor}$ & Upper: $2^{\left\lfloor\frac{n+1}{2}\right\rfloor}$ & Upper: $2^{\left\lfloor\frac{n+1}{2}\right\rfloor}$ \\
\hline
$(0,2)$ &
Tight: $n-2$~\cite{T-ITYM-Germany-South} & Tight: $n$~\cite{T-ITYM-Germany-South} & Tight: $n-1$~\cite{T-ITYM-Germany-South} \\
\hline
$(1,0)$ &
Tight: $n$ & Tight: $n$ & Tight: $n$ \\
\hline
$(2,1)$ &
Tight: $n$~\cite{T-ITYM-Germany-South} & Tight: $n-1$~\cite{T-ITYM-Germany-South} & Tight: $n-1$~\cite{T-ITYM-Germany-South} \\
\hline
\multirow{2}{*}{$(0,1)$} &
Lower: $\left\lfloor\frac{n-1}{2}\right\rfloor$ & Lower: $\left\lfloor\frac{n-1}{2}\right\rfloor$ & Lower: $\left\lfloor\frac{n}{2}\right\rfloor$ \\
& Upper: $n-1$ & Upper: $n$ & Upper: $n-1$ \\
\hline
\multirow{2}{*}{$(1,2)$} &
Lower: $\left\lfloor\frac{n}{2}\right\rfloor$ & Lower: $\left\lfloor\frac{n-1}{2}\right\rfloor$ & Lower: $\left\lfloor\frac{n-2}{2}\right\rfloor$ \\
& Upper: $n$ & Upper: $n$ & Upper: $n-1$ \\
\hline
\multirow{2}{*}{$(2,0)$} &
Lower: $\left\lfloor\frac{n}{2}\right\rfloor$ & Lower: $\left\lfloor\frac{n}{2}\right\rfloor$ & Lower: $\left\lfloor\frac{n}{2}\right\rfloor$ \\
& Upper: $n$ & Upper: $n$ & Upper: $n-1$ \\
\hline
\end{tabular}
\end{table}}

As we already mentioned, we know tight examples for infinitely many, but far not all, $n \equiv 2 \pmod 3$ for $b-a\equiv 1 \pmod 3$. The displayed universal lower bounds for these $(a,b)$ are based on the elementary examples where every element of $[n]$ belongs either to precisely one set or to all sets except one. 

When $a=b$ a classical block-type construction of size $2^{\left\lfloor\frac{n-m}{k}\right\rfloor}$ is
\begin{itemize}
        \item If $n\leq k$: $\mathcal{B}(m,k,n):=\{F\}$, where $F=[m]$.
        \item If $n>k$: for $N=\left\lfloor\frac{n-m}{k}\right\rfloor$ set
        \[\mathcal{B}(m,k,n):=
\Biggl\{[m]\cup \left(\bigcup_{i\in F}\{ki+m,\ ki+m-1,\ldots,\ ki+m-(k-1)\}\right)\ |\ F\subseteq\mathscr{P}(N)\Biggl\}\text{.}\]
\end{itemize}

The sets are formed as unions of $[m]$ and any number of blocks of size $k$. 

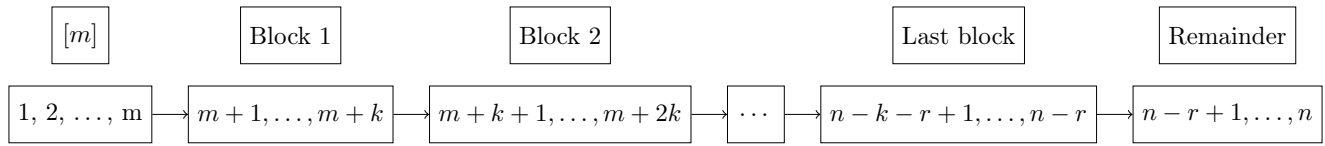
\begin{figure}[!h]
    \centering
\begin{center}
\resizebox{\textwidth}{!}{%
\begin{tikzpicture}[every node/.style={draw, minimum width=0.8cm, minimum height=0.8cm, align=center, font=\small}, node distance=0.5cm]
    \node[rectangle] (A) {1, 2, \dots, m};
    
    \node[rectangle, right=of A] (B1) {$m+1, \ldots, m+k$};
    \node[rectangle, right=of B1] (B2) {$m+k+1, \ldots, m+2k$};
    \node[right=of B2] (dots) {\dots};
    \node[rectangle, right=of dots] (Bn) {$n-k-r+1, \ldots, n-r$};
    
    \node[rectangle, right=of Bn] (R) {$n-r+1, \ldots, n$};
    
    \node[above=0.3cm of A] {$[m]$};
    \node[above=0.3cm of B1] {Block $1$};
    \node[above=0.3cm of B2] {Block $2$};
    \node[above=0.3cm of Bn] {Last block};
    \node[above=0.3cm of R] {Remainder};

    \draw[->] (A) -- (B1);
    \draw[->] (B1) -- (B2);
    \draw[->] (B2) -- (dots);
    \draw[->] (dots) -- (Bn);
    \draw[->] (Bn) -- (R);
\end{tikzpicture}
}
\label{fig:type-1-optimal}
\end{center}
\end{figure}

\vspace{-2em}

This construction is not maximal for moduli greater than $2$. Frankl and Odlyzko \cite{FRANKL-Odlyzko} gave a construction of size $(8k)^{\left\lfloor\frac{n}{4k}\right\rfloor}$ for $(0,0)$-town$\Mod k$ families for arbitrary $k\geq 2$, relying on Hadamard matrices. Furthermore, taking the unions of every set in the latter construction with a disjoint set of size $m$ yields an $(m,m)$-town of size $(8k)^{\left\lfloor\frac{n-m}{4k}\right\rfloor}$, also asymptotically larger than the block construction. Similar constructions exist for other generalizations of Eventown \cite{FRANKL-Uniform-Eventown-problems}.

Considering $(0,0)$-town$\Mod 3$, the Frankl-Odlyzko construction is a family of size $24^{\left\lfloor\frac{n}{12}\right\rfloor}$. A computer verification confirms that $24$ is the maximal size for $n=12$. See \cite{GitViz} for computed extremal sizes for all cases $(a,b)$ mod $k$, $k\in\{3,4,5,6\}$, up to $n=10$ (and up to $n=13$ for $(0,0)$-town$\Mod 3$).

\smallskip

Regarding the Ray-Chaudhuri--Wilson bound, we shall prove that there are infinitely many quadruples $(a,b,p,n)$ such that $|\mathcal{F}|\leq n-1$ in a sense where all of $a,b,n$ can vary for a fixed prime $p$.

\begin{theorem}
    \label{thm:infinitely-many-bound-less-than-n}
    Let $p$ be a prime and $a,b$ be integers with $0 \leq a, b \leq p-1$ and $p$ does not divide $b$ and $a-b$. There are infinitely many $n$ for which any $(a,b)$-town $\Mod p$ family $\mathcal{F}\subseteq \mathscr{P}(n)$ satisfies $|\mathcal{F}| \leq n-1$.  
\end{theorem}

We list notation used throughout the paper. For integers $x$ and $y$ write $y \mid x$ to indicate that $y$ divides $x$, otherwise write $y\nmid x$. For a set $A \in \mathscr{P}(n)$ denote its complement by $A^c = [n]\setminus A$. Denote by $\langle u, v \rangle := \sum_{i=1}^m u_i v_i$ the scalar product of $u = (u_1, \dots, u_m)$ and $v = (v_1, \dots, v_m)$. For a set $X$ and subset $A \subseteq X$ the indicator function $\mathbf{1}_A$ of $A$ on $X$ is such that $\mathbf{1}_A(x) = 1$ if $x\in A$ and $\mathbf{1}_A(x)=0$ otherwise. For a subset $A \in \mathscr{P}(n)$ denote by $ \chi_A := (\mathbf{1}_A(1), \mathbf{1}_A(2), \ldots, \mathbf{1}_A(n))$ its characteristic vector. 

\smallskip

It is sometimes useful to augment a characteristic vector by appending a scalar \(\alpha\) as its last coordinate. This construction can be used to encode additional information, independent of the vector's corresponding set. The idea is adapted from \cite{T-ITYM-Germany-South}.

\begin{defn}
    \label{defn:alpha-characteristic-vector}
    Let $\mathbb{F}$ be a field and $\alpha \in \mathbb{F}$ be a scalar. For a set $A \in \mathscr{P}(n)$ we denote by $\chi_A^\alpha$ the \emph{\(\alpha\)-characteristic vector} of $A$, where
    \[
        \chi_A^\alpha := \left(\mathbf{1}_A(1), \mathbf{1}_A(2), \ldots, \mathbf{1}_A(n), \alpha\right) \in \mathbb{F}^{n+1}.
    \]
\end{defn}

\section{Main proofs}

\subsection{Uniform upper bound for $a=b$} Recall that in a finite-dimensional vector space $V$, equipped with a symmetric bilinear form $\langle \cdot, \cdot \rangle$, a subspace $U\subseteq V$ is \emph{totally isotropic} if $\langle u_1, u_2 \rangle = 0$ for any (not necessarily distinct) $u_1, u_2 \in U$. It is well known that if the form is non-degenerate over $V$ and $U$ is totally isotropic, then $\dim U \leq \left\lfloor\frac{1}{2}\dim V\right\rfloor$. We also need (see e.g. \cite{Odlyzko1981}) that any $k$-dimensional vector space has at most $2^k$ binary vectors, i.e. such that each of their entries is $0$ or $1$.

\begin{proof}
    [Proof of Proposition \ref{prop:type-1-2-n-2-bound-intro}]
    Let $\alpha \in \mathbb{F}_{p^2}$ be such that $\alpha^2+m=0$. Denote by $\chi_i\in\mathbb{F}_{p^2}^{n+1}$ the $\alpha$-characteristic vector of $F_i\in\mathcal{F}$ and observe that $$\langle\chi_i,\chi_j\rangle=\alpha^2 + m = 0$$ for all $i,j$. If $U$ is the span of the $\chi_i$-s, then as $U$ is totally isotropic, we deduce $\dim U \leq \left \lfloor \frac{n+1}{2} \right\rfloor$. Finally, as the first $n$ entries of the $\alpha$-characteristic vectors are each $0$ or $1$, we conclude $|\mathcal{F}|\leq2^{\left\lfloor\frac{n+1}{2}\right\rfloor}$.
\end{proof}

\subsection{Upper bounds for $a\neq b$} A direct application of the Modular RW Theorem (Theorem~\ref{thm:modular-rw-theorem}) shows that the maximal size of an $(a,b)$-town$\Mod p$ family is bounded above by $n$ when $a\neq b$. Let us present an alternative proof using $\alpha$-characteristic vectors, whose ideas will be useful for obtaining the stronger results.

\begin{lemma}
\label{lemma:a-neq-b-partial-bound}
    Let $\mathcal{F}\subseteq\mathscr{P}(n)$ be an $(a,b)$-town$\Mod p$, where $p$ does not divide $a-b$. Then $|\mathcal{F}|\leq n$.
\end{lemma}

\begin{proof}
    Let $\alpha$ be an element of $\mathbb{F}_{p^2}$ such that $\alpha^2+b=0$. Let $\mathcal{F}=(F_1,\ldots,F_{|\mathcal{F}|})$. Denote by $\chi_i\in\mathbb{F}_{p^2}^{n+1}$ the $\alpha$-characteristic vector of $F_i\in\mathcal{F}$ for each $i=1,\ldots,|\mathcal{F}|$, and observe that
    \[
        \langle\chi_i,\chi_j\rangle = |F_i\cap F_j|+\alpha^2 = 
        \begin{cases}
            a-b, & \text{if } i=j, \\
            0, & \text{if } i\neq j.
        \end{cases}
    \]
    Let $\lambda_1,\ldots,\lambda_{|\mathcal{F}|}$ be scalars such that $\lambda_1\chi_1+\ldots+\lambda_{|\mathcal{F}|}\chi_{|\mathcal{F}|}=0\text{.}$ Taking the scalar product with $\chi_i$, we obtain $(a-b)\lambda_i=0$ and hence $\lambda_i = 0$ for all $i$, i.e. the vectors are linearly independent. Moreover, each of them is orthogonal to $v = (1,1,\ldots,1,-a\alpha^{-1})$ for $b\not\equiv 0 \Mod p$ and to $e = (0,0,\ldots,0,1)$ for $b\equiv 0 \Mod p$, so they lie in a subspace $V \subseteq \mathbb{F}_{p^2}^{n+1}$ of dimension $n$. The result follows.
\end{proof}

We are now ready to proceed to the improvements of the upper bound to $n-1$.

\begin{proposition}
\label{prop:summary-n-1}
    Let $a,b,n$ be non-negative integers and $p$ be a prime. Suppose at least one of the following holds$:$
    \begin{enumerate}[label=\emph{(\roman*)}]
        \item $p \nmid a, \ b, \ a-b, \ n, \ a^2-nb-a+b;$
        \item $p\mid a$ and $p\nmid b, \ n-1$.
    \end{enumerate}
    Then for any $(a,b)$-town$\Mod p$ family $\mathcal{F} \subseteq \mathscr{P}(n)$ we have $|\mathcal{F}| \leq n-1$.
\end{proposition}

Note that Theorem $\ref{thm:infinitely-many-bound-less-than-n}$ follows for $p\geq 3$ from this proposition, as for any $a$ and $b$ with $p \nmid b, a-b$ there are at most two forbidden congruence classes for $n \Mod p$, hence at least one attainable.  

\begin{proof}
    In both cases the conditions of Lemma~\ref{lemma:a-neq-b-partial-bound} hold, so throughout we shall use all notation and computations from its proof. We always have $|\mathcal{F}| \leq n$.

    Suppose firstly that $p\nmid a, \ b, \ a-b, \ n, \ a^2-nb-a+b$. Assume for contradiction that $|\mathcal{F}|=n$. Then $\chi_1,\ldots,\chi_n$ form a basis of $V$. Consider $u=(1,1,\ldots,1,na^{-1}\alpha)$ and observe that $u\in V$ since $p\nmid b$ and $u$ is orthogonal to $v$. Then there exist $\lambda_1,\ldots,\lambda_n$, at least one being nonzero, such that
    $u=\lambda_1\chi_1+\ldots+\lambda_n\chi_n$.
    Taking the scalar product of both sides with $\chi_i$, we obtain for all $i=1,\ldots,n$
    \[a-na^{-1}b=\lambda_i(a-b)\Leftrightarrow\lambda_i=\frac{a-na^{-1}b}{a-b}=\frac{a^2-nb}{a(a-b)}\text{.}\]
    Denote the common value of $\lambda_i$ by $\Lambda$. If $p\mid a^2-nb$, then $\lambda_i = \Lambda = 0$ for all $i$, contradiction and we are done. Suppose $p\nmid a^2-nb$. Now that $\displaystyle u=\Lambda\sum_{i=1}^n\chi_i$, we derive another expression for $\Lambda$, this time by taking the scalar product of $u$ with itself. This gives
    \[n-n^2a^{-2}b  =\langle u,u\rangle=\Lambda^2\left\langle\sum_{i=1}^n \chi_i, \sum_{j=1}^n \chi_j\right\rangle=\Lambda^2n(a-b)\text{.}\]
    Therefore, as $p\mid n$,
    \[ \left(\frac{a^2-nb}{a(a-b)}\right)^2 = \Lambda^2 = \frac{1 - na^{-2}b}{a-b} = \frac{a^2-nb}{a^2(a-b)}.\]
    Cancelling common non-zero terms on both sides now leads to $a^2 - nb = a-b$ in $\mathbb{F}_{p^2}^{n+1}$. This contradicts the assumption $p\nmid a^2-nb-a+b$ and completes the proof in this case.

    \smallskip

    Now suppose $p\mid a$ and $p \nmid b, n-1$. Assume for contradiction that $|\mathcal{F}|=n$. Then $\chi_1,\ldots,\chi_n$ form a basis of $V$. Consider $e = (0,0,\ldots,0,1)$ and observe that $e\in V$ since $v=(1,1,\ldots,1,0)$ when $p\mid a$. Then there exist $\mu_1,\ldots,\mu_n$, at least one being nonzero, such that
    $e=\mu_1\chi_1+\ldots+\mu_n\chi_n$.
    Taking the scalar product of both sides with $\chi_i$, we obtain $\alpha = \mu_i(a-b) = -\mu_ib$, i.e. $\mu_i=-\alpha b^{-1}$. Denote the common value of $\mu_i$ by $M$. Now that $\displaystyle e = M\sum_{i=1}^n\chi_i$, we derive another expression for $M$, this time by taking the scalar product of $e$ with itself -- or equivalently, by comparing the last coordinate in $\displaystyle e = M\sum_{i=1}^n\chi_i$. We obtain $1 = Mn\alpha$. Substituting $M=-\alpha b^{-1}$ leads to $1 = -n\alpha^2b^{-1} = n$ in $\mathbb{F}_{p^2}^{n+1}$. This contradicts the assumption $p\nmid n-1$ and completes the proof.
\end{proof}

We now check how the $(a,b)$-town property is changed when taking complements.

\begin{defn}
\label{defn:substitution}
    For a family $\mathcal{F}= \left\{F_1,\ldots,F_m\right\}\subseteq\mathscr{P}(n)$ the \emph{substitution family} is $\mathcal{F}_\xi:=\left\{F_1^{c},\ldots,F_m^{c}\right\}$. 
\end{defn}

Since $|\mathcal{F}| = |\mathcal{F}_\xi|$, any bound for $\mathcal{F}_\xi$ is also a bound for $\mathcal{F}$. In what follows it is important that the difference $a-b$ is the same in the original and in the substitution family.

\begin{lemma}
\label{lemma:substitution-families-nums}
    Let $\mathcal{F}\subseteq\mathscr{P}(n)$ be an $(a,b)$-town$\Mod k$. Then $\mathcal{F}_\xi$ is an $(n-a, n-2a+b)$-town$\Mod k$.
\end{lemma}

\begin{proof}
    Let $\mathcal{F}=\left\{F_1,\ldots,F_m\right\}$ and $\mathcal{F}_\xi=\left\{F_1^{c},\ldots,F_m^{c}\right\}$. We have $\left|F_i^{c}\right| = n - \left|F_i\right| \equiv n-a\Mod k$ and
    \[\left|F_i^{c}\cap F_j^{c}\right|\equiv n-\left|F_i\right|-\left|F_j\right|+\left|F_i\cap F_j\right|\equiv n-2a+b\Mod k\text{, }\ i\neq j\text{.}\qedhere\]
\end{proof}

Combining Proposition~\ref{prop:summary-n-1} and Lemma~\ref{lemma:substitution-families-nums} immediately gives even more tuples $(a,b,p,n)$, where the upper bound is less than $n$.

\begin{proposition}
\label{prop:n-1-subs}
    Let $a,b,n$ be non-negative integers and $p$ be a prime. Suppose at least one of the following holds$:$
    \begin{enumerate}[label=\emph{(\roman*)}]
        \item $p \nmid n-a, \ n-2a+b, \ a-b, \ n, \ a^2-nb-a+b;$
        \item $p\mid n-a$ and $p\nmid n-2a+b, \ n-1$.
    \end{enumerate}
    Then for any $(a,b)$-town$\Mod k$ family $\mathcal{F} \subseteq \mathscr{P}(n)$ we have $|\mathcal{F}| \leq n-1$.
\end{proposition}

It is interesting that the conditions $p\nmid a-b$ and $p\nmid a^2 - nb - a + b$ do not change when taking complements, and it would be nice to understand a combinatorial reason behind the second one.

\smallskip

Now we deduce all new improvements for modulo $3$.

\begin{proof}
    [Proof of Proposition~\ref{prop:mod-3-n-1}] Suppose $n\equiv 0 \Mod 3$, then $a=0, b=1$ follows from the second part of Proposition~\ref{prop:summary-n-1}. Now suppose $n\equiv 2 \Mod 3$. Here $a=0, b=1$ follows from the second part of Proposition~\ref{prop:summary-n-1}, while $a=1,b=2$ follows from the first part of Proposition~\ref{prop:summary-n-1} and $a=2, b=0$ follows from the second part of Proposition~\ref{prop:n-1-subs}.
\end{proof}

Apart from the classical Oddtown $a=1$, $b=0$, we can deduce tight bounds when $a=2$, $b=1$ for $n$ and $p$ only related by a congruence condition. In fact, only $n\not\equiv 0 \Mod p$ is currently out of reach, so it would be great if a version of Proposition \ref{prop:summary-n-1} with $p\mid n$ can be derived in order to cover this.

\begin{corollary}
    \label{corollary:(1,2)-n-p-congruence-related}
    Let $\mathcal{F}\subseteq\mathscr{P}(n)$ be a $(2,1)$-town$\Mod p$. If $n\equiv 3\Mod p$, then $|\mathcal{F}|\leq n$, and if $n\not\equiv 0,3\Mod p$, then $|\mathcal{F}|\leq n-1$. In all cases the bound is tight.
\end{corollary}

\begin{proof}
    When $a=2$, $b=1$, $n\equiv 3 \Mod p$ the substitution family $\mathcal{F}_\xi$ is a $(1,0)$-town$\Mod p$, since $n-a \equiv 3-2 \equiv 1$ and $n-2a+b \equiv 3-4+1 \equiv 0$. By the classical Oddtown problem for a field of characteristic $p$, we have that $n$ is a tight upper bound for the $(1,0)$-town $\mathcal{F}_\xi$, hence also for $\mathcal{F}$. 

    When $a=2$, $b=1$, $n\not\equiv 0,3 \Mod p$ the bound follows by Proposition~\ref{prop:summary-n-1} for $a=2,b=1$ and is attained by the family $\{\{1,2\},\{1,3\},\ldots,\{1,n\}\}$.
\end{proof}

We conclude with a short discussion of hypotheses based on modest numerical evidence from \cite{GitViz}. Here the integer $k$ need not be prime.

\begin{conjecture}
    For $i=1,2$ let $\mathcal{F}_i\subseteq\mathscr{P}(n)$ be a $(m_i,m_i)$-town$\Mod k$ family of maximum size, where $0 \leq m_1, m_2 \leq k-1$. If $m_1 < m_2$, then $|\mathcal{F}_1|\geq |\mathcal{F}_2|\text{.}$
\end{conjecture}

We already know that the linear upper bound $n$ holds when the $(a,b)$-town satisfies $a\not\equiv b \pmod p$ for some prime divisor $p$ of $k$ -- this holds e.g. when $a\neq b$ and $k$ is squarefree. What remains is where $a\equiv b\Mod p$ for every prime divisor $p\mid k$, among which we expect that $a\equiv b \Mod k$ is the sole exceptional case of exponential size and the rest are linear.

\begin{conjecture}
    Let $\mathcal{F}\subseteq\mathscr{P}(n)$ be an $(a,b)$-town$\Mod k$, where $a\not\equiv b \Mod k$. Then $|\mathcal{F}|\leq n$.
\end{conjecture}

\nocite{*}
\bibliographystyle{plain}
\bibliography{refs}

 \end{document}